\documentclass[3p]{elsarticle}

\usepackage{amsmath}
\usepackage{amssymb}
\usepackage{hyperref}
\usepackage{graphicx}

\usepackage{subcaption}

\usepackage{siunitx}
 
\newcommand{\itemEq}[1]{%
         \begingroup%
         \setlength{\abovedisplayskip}{0pt}%
         \setlength{\belowdisplayskip}{0pt}%
         \parbox[c]{\linewidth}{\begin{flalign}#1&&\end{flalign}}%
         \endgroup}
\newcommand{\defint}[4]{\int_{#1}^{#2}{#3} d{#4}}
\newcommand{\w}{\omega}

\newcommand{\expn}[1]{e^{#1}}

\newcommand{\Ptilde}{\tilde{P}}

\usepackage{longtable} 
\usepackage{mathtools}

\newcommand{\appropto}{\mathrel{\vcenter{
  \offinterlineskip\halign{\hfil$##$\cr
    \propto\cr\noalign{\kern2pt}\sim\cr\noalign{\kern-2pt}}}}}
\usepackage[utf8]{inputenc}

\usepackage{siunitx,amsmath}

\usepackage{xcolor}
\newcommand{\clearpageSamiya}{ }

\journal{Applied Mathematics and Computation}

\begin{document}

\begin{frontmatter}

\title{Analytic Methods for Differential Algebraic Equations}

% \title{Methods for the Analytic Study of Sys Differential Algebraic Equations}

\author{Samiya A Alkhairy \footnote{Correspondence to samiya@mit.edu, samiya@alum.mit.edu}}
\address{Earth, Atmospheric, and Planetary Sciences Department \\ Massachusetts Institute of Technology \\ 77 Mass Ave Office 54-317, Cambridge, MA 02139}

\begin{abstract}

We introduce methods for deriving analytic solutions from differential-algebraic systems of equations (DAEs), as well as methods for deriving governing equations for analytic characterization which is currently limited to very small systems as it is carried out by hand. Analytic solutions to the system and analytic characterization through governing equations provide insights into the behaviors of DAEs as well as the parametric regions of operation for each potential behavior.  DAEs are a mixture of algebraic, ordinary differential, and/or partial differential equations which are multivariate constitutive equations that arise from fundamental laws of a discipline which relate multiple dependent variables. For each system (DAEs), and choice of dependent variable, there is a corresponding governing equation which is univariate ODE or PDE that is typically higher order than the constitutive equations of the system. We first introduce a direct formulation for representing systems of linear differential-algebraic scalar constitutive equations that are homogeneous or nonhomogeneous with constant or nonconstant coefficients. Unlike state space formulations, our formulation follows very directly from the system of constitutive equations without the need for introducing state variables or singular matrices. Using this formulation for the system of constitutive equations (DAEs), we develop methods for deriving analytic expressions for the full solution (complementary and particular) for all dependent variables of systems that consist of constant coefficient ordinary-DAEs and special cases of partial-DAEs. We also develop methods for deriving the governing equation for a chosen dependent variable for the constant coefficient ordinary-DAEs and partial-DAEs as well as special cases of variable coefficient DAEs. The methods can be automated with symbolic coding environments thereby allowing for dealing with systems of any size while retaining analytic nature. This is relevant for interpretable modeling, analytic characterization and estimation, and engineering design in which the objective is to tune parameter values to achieve specific behavior. Such insights cannot directly be obtained using numerical simulations.

\end{abstract}


\end{frontmatter}

% \linenumbers

\clearpageSamiya

% \tableofcontents

% \clearpageSamiya

\section{Introduction}

% \cite{Feynman1963118,Dirac1953888}

\subsection{Constitutive Equations and Governing Equations}
\label{s:constitutiveGoverningMotivation}

Many systems in engineering, science, and social science are described by fundamental principles and laws. These laws are expressed as simultaneous equations in terms of multiple variables dependent on time and / or space. We refer to these equations as constitutive equations. These multivariate equations relate the various dependent variables of a system and, together, fully describe said system. By the nature of the laws resulting in constitutive equations, the problem is typically perfectly determinable: There are as many equations as there are unknowns (dependent variables) and all equations are independent of one another. An equation is not a combination of the other equations but rather arises from a distinct law. The set of constitutive equations that describe a system is referred to as the system of constitutive equations. The systems we are interested in consist of ordinary differential, partial differential and/or algebraic equations. The differential equations in these systems may be constant or nonconstant coefficient, homogeneous or nonhomogeneous, and of any order. Some such systems are linear, and others are nonlinear (potentially linearized). The system may have one or multiple independent variables.  Constitutive systems of differential-algebraic equations arise in a variety of problems such as multibody mechanics \cite{Schutte, Pappalardo}, membrane distillation \cite{membraneDistillation1, membraneDistillation2}, chemical kinetics \cite{Beigler, Balawi} and circuit modified node analysis and other network equation methods \cite{originalMNA, MNA, MNA2, NilpotencyIndex}. Therefore, it is important to develop analytical (and numerical) methods for studying these systems. Domain experts may be interested in: determining solutions via analytic methods or numerical simulation, interpretations and characterization through analytic study such as determining possible behaviors or parametric regions of operation, studying control aspects of these systems, or parameter estimation problems. 

We are particularly interested in developing analytic methods for: (a) obtaining analytic expressions for solutions to the system, as well as (b)  deriving governing equations - which may be thought of as decoupling or reducing to a single dependent variable. 

Reduction of a system of constitutive equations into a governing equation for a given dependent variable produces a single univariate ODE or PDE (uni-ODE or uni-PDE) that is of equal or higher order compared to the constitutive equations \footnote{A governing equation is specific to a system of constitutive equations and a single chosen dependent variable. In the case of homogeneous constant coefficient systems, the governing equation is the same for all dependent variables}.  Examples of governing equations derived from constitutive equations include the wave equation, Helmholtz equation, and the governing equation describing the displacement of a particular mass in a vibrating system of masses.

The structure of the governing equation is an extremely powerful analytic tool for studying a system systems. Consider, for instance, the classification of governing equations that are partial differential equations into elliptic/parabolic/hyperbolic. A system-specific example for the benefit of studying the governing equation is the well-known Webster horn equation which is derived from the constitutive equations in order to study a horn's properties which is of relevance for a variety of systems. Another particularly simple example is that combining equations to arrive at the governing equation for ventilator-respiratory system is key to developing a noninvasive estimation method for determining patient respiratory system parameters and flow variables from measurements at the ventilator \cite{philipsPaper}. These examples motivate the need for developing a systematic method for deriving the governing equation from the system of constitutive equations. Despite this, the derivation of the governing equation from differential-algebraic equations is typically done by hand and is therefore only derived if the number of constitutive equations is very small (generally consisting of 2-3 dependent variables).

% The governing equation is usually only derived if the number of constitutive equations is very small (generally consisting of two dependent variables). This is because the constitutive equations are usually combined by hand to form the governing equation which is then used for analytic characterizing and for solving for analytic solutions.

The form of the analytic solutions to constitutive systems of differential-algebraic equations (sys-DAEs) is also key to studying these systems. Current methods to obtain analytic solutions exist for certain classes of sys-DAEs. These require first formulating the system using the general state-space formulation \footnote{Alternatively, the solution for a single dependent variable of the system is determined from the governing equation after it is derived from the constitutive equations by hand. In this case, single equation (uni-ODE) methods such as those of \ref{s:differentialOperator} are used.}. The formulation does not retain the natural structure of the system but requires introducing state variables. Current analytic solution methods are based on the state-space formulation, and a survey of these solution methods is provided in \ref{s:literatureSurvey}. We do not build on the state-space formulation and instead propose an alternative representation for the analytic study of systems, but have included this appendix as the state-space formulation is widely used. Abbreviations we use throughout this paper for types of constitutive and governing equations are in Section \ref{s:abbreviations}.

In this paper, we present a direct formulation for the system of constitutive equations in its natural variables. We then develop methods to facilitate analytical investigation of both the constitutive equations and the governing equation. Derivation of governing equations and analytic expressions for solutions generally enable us to understand system behavior in ways inaccessible by numerical simulations. The formulation and methods are applicable for systems of any size, and are programmable into a symbolic computing environment. In this paper, we restrict ourselves to constitutive systems of smooth linear ordinary/partial differential-algebraic equations (sys-LDAEs). We only concern ourselves with scalar differential operators - as opposed to operations such as div, grad, and curl. This paper is not concerned with studying these systems from purely mathematical and systems theoretic perspectives. The formulation and methods we present in this paper are for differential-algebraic systems of equations and do not require any modifications for the subset of these systems that only contain differential equations.

% Problems in the time-domain are mostly initial value problems. More generally, we may have initial, final, boundary, and mixed value problems, especially when spatial domains are involved. 

\clearpageSamiya

\subsection{Abbreviations}
\label{s:abbreviations}

We use the following abbreviations for univariate equations (uni) and multivariate systems of equations (sys) in this paper. All equations may be homogeneous or nonhomogeneous. As mentioned previously, we are only concerned with linear (L) equations and cases with scalar differential operators and do not consider vector differential operators. All systems we consider consist of independent equations which are equal in number to the number of dependent variables.  

\begin{table}[htbp]
    \centering
    \begin{tabular}{|l|l|}
    \hline
    \textbf{Abbreviation} & \textbf{Definition} \\ \hline
    uni-LAE & Algebraic equation \\ \hline
    uni-LODE     & Ordinary differential equation \\ \hline
    uni-LPDE     & Partial differential equation \\ \hline
    sys-LAE     & System of  algebraic equations \\ \hline
    sys-LODEs     & System of  ordinary differential equations \\ \hline
    sys-LPDEs     & System of  partial differential equations \\ \hline
    sys-LODAEs     & System of  ordinary differential-algebraic equations \\ \hline
    sys-LPDAEs     & System of  partial differential-algebraic equations \\ \hline
    sys-LDAEs    & sys-LODAEs and sys-LPDAEs \\ \hline
    cc (prefix, not needed for LAEs) & Constant coefficient - e.g. uni-ccLODEs \\ \hline
    vc (prefix, does not apply to LAEs) & Variable coefficient - e.g. sys-vcLODAEs \\ \hline

    \end{tabular}
    \caption{Abbreviations. Note that in the literature the abbreviation DAE usually refers to sys-ODAEs.}
    \label{tab:abbreviations}
\end{table}

\clearpageSamiya

\subsection{Objectives}
\label{s:objectives}

In Section \ref{s:constitutiveGoverningMotivation}, we motivated the need for methods for analytically solving a system of constitutive equations as well as methods for deriving governing equations from systems of constitutive equations. More explicitly, our objectives are to develop the following:

\begin{enumerate}
    \item Objective I: Formulation for systems of constitutive equations that are sys-LDAEs.
    \item Objective II: Systematic methods for deriving governing equations (uni-LODE / uni-LPDE) from sys-LDAEs for a chosen dependent variable.
    \item Objective III: Systematic methods for deriving analytic expressions for full solutions of dependent variables from sys-LDAEs.
\end{enumerate}

The constant coefficient class of systems is the natural starting point to develop any analytic methods as evidenced by the literature. Hence, over the course of this paper, we focus on developing the  formulation (objective I) and simple analytic methods (objectives II, III) for homogeneous or nonhomogeneous sys-ccLODAEs that include multivariate algebraic equations and ODEs of any order. We then discuss methods and limitations for other types of sys-LDAEs (sys-ccLPDAEs, sys-vcLODAEs, and sys-vcLPDAEs).

The constitutive equations can unambiguously and directly be represented using the proposed formulation - as opposed to having to define state variables.  The formulation further directly admits algebriac equations and allows for partial differential equations. We use the formulation to develop methods for deriving the governing equation and obtaining analytic solutions to the system in its native domains. The methods are suitable for systems of any size and  can be programmed into an symbolic system such as Maple and Mathematica using existing linear algebraic algorithms.  The simplicity and the ability to encode the proposed analytic methods in symbolic toolboxes enables them to be used by domain experts who encounter sys-LDAEs in their systems (for example in certain distillation equations or linear circuit modified nodal analysis) rather than being limited in use to those specialized in mathematics or systems theory.

 %Existing analytic methods for obtaining solutions are quite simple for first order sys-ccLODEs (see Section \ref{s:classicalFormulations}) but are generally otherwise complicated or specific to certain problems.

% Our formulation is not restricted to initial value problems. 

\clearpageSamiya

\subsection{Paper Organization and Examples}

In developing the formulation and methods towards the objectives in Section \ref{s:objectives}, we build on concepts from linear algebra and differential operator theory. Due to this, and the fact that the methods introduced in this paper are intended to be used in a variety of fields by engineers, scientists, and social scientists, we review differential operator theory concepts we use in \ref{s:differentialOperator}. We refer to equations from the appendix throughout the paper and the reader is encouraged to start with \ref{s:differentialOperator} if unfamiliar with these concepts. Also included in \ref{s:differentialOperator}, are some modifications we included in differential operator theory that allow our methods to be more direct and general than they otherwise would be. 

We introduce our formulation and methods starting with sys-ccLODAEs in Sections \ref{s:formulate}, \ref{s:governingMethod}, and \ref{s:solutionMethod}. We then introduce the formulation, methods, and limitations for sys-ccLPDAEs in Section \ref{s:ccLPDAEs}, and finally for sys-vcLODAEs and sys-vcLPDAEs in Section \ref{s:vc}. Throughout this paper, we provide many examples so that the formulation and methods can be better understood and applied. The examples include transmission lines and vibrating mass systems.

\clearpageSamiya

\section{Formulating Constitutive sys-ccLODAEs as a Linear System of Equations}
\label{s:formulate}

We introduce a formulation for the sys-LDAE system of constitutive equations (objective I), and focus here, in particular, on the sys-ccLODAEs. The formulation is specifically constructed for developing the analytic methods for objectives II and III (deriving the governing equation and obtaining analytic solutions) rather than, for instance, numerical simulations, study of properties from controls perspectives, or determining eigenvalues and eigenvectors for certain classes of equations.

\subsection{Proposed Formulation for Objectives II and III}

Consider a sys-LAE,
\begin{equation}
    A x = b
\end{equation}

\begin{equation}
\begin{bmatrix}
l_{11} & \cdots & l_{1n}\\
\vdots &  \ddots & \vdots \\
l_{n1} & \cdots & l_{nn}
\end{bmatrix} 
\begin{bmatrix} x_1 \\ \vdots \\ x_n \end{bmatrix} 
= 
\begin{bmatrix} b_1\\ \vdots \\ b_n \end{bmatrix} 
\label{eq:sysLAE}
\end{equation}

where we are only concerned with the case of a nonsingular square matrix $A$ which is the category relevant for constitutive equations.

Our formulation for the general class of sys-LDAEs involves expressing the system of constitutive equations as if it were a sys-LAE with each row corresponding to a constitutive equation. For sys-ccLODAEs, this is,

\begin{equation}
    M(D) x(t) = f(t)
    \label{eq:ccLODAEsForm}
\end{equation}

or, more explicitly,

\begin{equation}
\begin{bmatrix}
P_{11}(D) & \cdots & P_{1n}(D)\\
\vdots &  \ddots & \vdots \\
P_{n1}(D) & \cdots & P_{nn}(D)
\end{bmatrix} 
\begin{bmatrix} x_1(t) \\ \vdots \\ x_n(t) \end{bmatrix} 
= 
\begin{bmatrix} f_1(t) \\ \vdots \\ f_n(t) \end{bmatrix} 
    \label{eq:ccLODAEsFormDetails}
\end{equation}

with $t$ as the independent variable. $M(D)$ is a matrix where the entries are polynomials of differential operator $D \mapsto \frac{d}{dt}$; $x(t)$ is a vector of time-dependent variables, $x_i(t)$; and $f(t)$ is a vector of forcing terms $f_i(t)$.

By nature of the set of the inherent laws / constitutive equations for a typical system of interest, the matrix of operators will always be square and the rows are all independent. Note that AEs are included in this form by assigning all polynomials in the corresponding row to be of zeroth order. We also note that many integral and integro-differential constitutive equations may be formulated as equivalent differential constitutive equations by differentiation (though this may lead to some loss of information). Therefore, our formulation as a system of linear equations - objective I, and methods for objectives II and III also extend to systems that include constitutive integro-differential equations that can be formulated accordingly (by multiplication of the corresponding row with some $D^k$). Any $\frac{df_i}{dt}$ terms are absorbed into $f$ and there is no need for defining ant new intermediate variables.

% Formulating sys-ccLODAEs as above enables us to, and is the first step towards, developing and applying methods for objectives II and III. The formulation is similar for the larger class of sys-LDAEs as discussed in sections \ref{s:ccLPDAEs} for sys-ccLPDAEs and \ref{s:vc} for sys-vcLODAEs and sys-vcLPDAEs. 

Our formulation is appropriate for, and is the first step for, deriving the governing equation (objective II) for constitutive sys-ccLODAEs, sys-ccLPDAEs and certain classes of sys-vcLODAEs and sys-vcLPDAEs; deriving an analytic expression for the particular and full solutions that does not require simulation (objective III) as we discuss for ccLODAEs and certain classes of ccLPDAEs; and analytically characterizing the system and obtaining physical insights. The methods are appropriate for the simple, systematic and automatable analytic study of systems of any size.

\subsection{Examples}

In this section, we exemplify our formulation for sys-ccLODAEs. We include the equivalent state-space formulation for these examples.

For our first example, consider the following set of constitutive equations,
\begin{equation}
    \begin{cases}
    \ddot{x}_1(t) - x_2(t) = f_1(t)\\
    \dot{x}_2(t) - x_3(t) = f_2(t)\\
    \dot{x}_3(t) + x_3(t) - x_1(t) + x_2(t) = f_3(t)
    \end{cases}
\end{equation}

Our formulation (Equation \ref{eq:ccLODAEsForm}) follows directly from the above and is,
\begin{equation}
\begin{bmatrix}
D^2 & -1 & \\
 & D & -1 \\
-1 & 1 & D + 1 
\end{bmatrix} 
\begin{bmatrix} x_1 \\ x_2 \\ x_3 \end{bmatrix} 
= 
\begin{bmatrix}
f_1 \\ f_2 \\ f_3
\end{bmatrix}
\end{equation}

Such a class of systems also has an easily obtainable corresponding state-space formulation of equation  \ref{eq:stateSpace1}. If the system contains an equation with a second order derivative in one of the dependent variables, then, for the state-space formulation, a new dependent variable is defined and the defining equation added to the system of equations. For the above example, this leads to,
\begin{equation}
\dot{\begin{bmatrix} x_1 \\ x_2 \\ x_3 \\ x_4 \end{bmatrix} }
= 
\begin{bmatrix}
&&&1 \\
 &  & 1 & \\
1 & -1 & -1 & \\ 
&1&&
\end{bmatrix} 
\begin{bmatrix} x_1 \\ x_2 \\ x_3 \\ x_4  \end{bmatrix} 
+ 
\begin{bmatrix}
0 \\ f_2 \\ f_3 \\ f_1
\end{bmatrix}
\end{equation}

For our second example, consider the following set of constitutive equations,
\begin{equation}
    \begin{cases}
    \ddot{x}_1(t) - x_2(t) + x_4(t) = f_1(t) \\
    -x_1(t) + \dot{x}_2(t) + \dot{x}_3(t) = f_2(t) \\
    \dot{x}_2(t) + x_2(t) + x_4(t) = f_3(t) \\
    x_1(t) - x_2(t) - x_3(t) + x_4(t) = f_4(t)
    \end{cases}
\end{equation}

Our formulation (Equation \ref{eq:ccLODAEsForm}) which follows directly from the above system is,
\begin{equation}
\begin{bmatrix}
D^2 & -1 & & 1\\
-1 & D & D & \\
 & D+1 &  & 1 \\
 1 & -1 & -1 & 1
\end{bmatrix} 
\begin{bmatrix} x_1 \\ x_2 \\ x_3 \\ x_4 \end{bmatrix} 
= 
\begin{bmatrix}
f_1 \\ f_2 \\ f_3 \\ f_4
\end{bmatrix}
\end{equation}

The formulation for the same set of equations using singular state-space (Equation \ref{eq:stateSpaceDAE}) is,

\begin{equation}
    \begin{bmatrix}
        1&&&& \\
        &&&&1 \\
        &1&1&& \\
        &1&&& \\
        &&&&
    \end{bmatrix}
    \dot{\begin{bmatrix} x_1 \\ x_2 \\ x_3 \\ x_4 \\ x_5 \end{bmatrix} }
    = 
    \begin{bmatrix}
        &&&&1 \\
        &1&&-1& \\
        1&&&& \\
        &-1&&-1& \\
        -1&1&1&-1&
    \end{bmatrix}
    \begin{bmatrix} x_1 \\ x_2 \\ x_3 \\ x_4 \\ x_5  \end{bmatrix} 
    + 
    \begin{bmatrix} 0 \\ f_1 \\ f_2 \\ f_3 \\ f_4  \end{bmatrix} 
\end{equation}

The matrix multiplying $\dot{x}(t)$ encodes the linear algebraic equation in the final row, and the constant coefficient ODE that contain derivatives in more than one dependent variable in the third row. 

The choice of formulation is the basis for developing methods for objectives II and III. Notice that the proposed formulation follows unambiguously from the constitutive equations. There is no introduction of intermediate variables (state variables) to study a system with higher order derivatives. The proposed formulation does not involve singular matrices which are traditionally needed for algebraic equations and equations with derivatives in more than one state variable. As discussed in Section \ref{s:ccLPDAEs}, the formulation also admits partial differential equations.

In \ref{s:literatureSurvey}, we review existing methods for objective III that are for the general state-space formulation. The rest of the paper is concerned with developing methods for objectives II and III based on our direct formulation.

\clearpageSamiya

\section{Methods for Deriving Governing uni-ccLODE from Constitutive sys-ccLODAEs} 
\label{s:governingMethod}

In this section, we develop methods for objective II - systematically deriving the higher (or equal) order governing uni-ccLODE for a particular dependent variable from a system of constitutive equations sys-ccLODAEs formulated as in Equation \ref{eq:ccLODAEsForm}.

\subsection{Method}
\label{s:governingMethodPart}

In order to describe our method for objective II, consider Gaussian elimination for sys-LAEs (Equation \ref{eq:sysLAE}). Gaussian elimination involves a series of row operations which leads to a final row,
\begin{equation}
    h(l) x_n = \sum_{i=1}^n g_i(l)b_i
    \label{eq:gaussElimResultLinAlg}
\end{equation}
where $l$ is a vector containing the entries $l_{ij}$ of matrix $A$. And $h$ and $g_i$ are functions of their arguments in a form dictated by Gaussian elimination row operations. Row permutation is occasionally needed during Gaussian elimination.

To derive higher order governing uni-LODEs/uni-LPDEs in a specified dependent variable, we use Gaussian elimination. Let $x_n(t)$ be the variable of interest for the governing equation. We formulate the system as in Equation \ref{eq:ccLODAEsForm} with the variable of interest as the final variable \footnote{Related to this is the fact that the higher order governing equation for homogeneous sys-ccLODAEs is the same for any chosen dependent variable. However, we explicitly state choosing the dependent variable to be the final variable in the vector $x$ here as it is necessary to specify for the nonhomogeneous case, for the case in which $M(D)$ of Equation \ref{eq:ccLODAEsForm} can be written as block diagonal, and for the variable coefficient case.}. 

We then carry out steps of Gaussian elimination on the system, which yields a parallel to Equation \ref{eq:gaussElimResultLinAlg},

\begin{equation}
    h(\Ptilde) x_n(t) = \sum_{i=1}^n g_i(\Ptilde) f_i(t)
    \label{eq:gaussccLODAEs}
\end{equation}
where $\Ptilde$ is a vector of the polynomial operators $P_{ij}(D)$ which are entries of matrix $M(D)$ of Equation \ref{eq:ccLODAEsForm}. $h$ and $g_i$ have forms arising from row operations. Note that we do not need to keep track of the order of operations of Gaussian elimination to reach Equation \ref{eq:gaussccLODAEs} for the case of sys-ccLODAEs and they may be considered commutative because the equations are constant coefficient. This arises because $Dax = aDx$. However, this will not generally be the case for sys-vcLODAEs and sys-vcLPDAEs as we discuss in Section \ref{s:vcGoverning}. 

We then manipulate this expression into the governing uni-ccLODE by multiplication of both sides by all denominators of the $n^{th}$ rows of the elimination matrices used in Gaussian elimination to obtain Equation \ref{eq:gaussccLODAEs} to yield the following form. This manipulation via multiplication is possible owing to the constant coefficients nature.

\begin{equation}
\mathcal{P}(D) x_n(t) = \sum_{i=1}^{n} \mathcal{Q}_i(D) f_i(t)
\end{equation}

Where $\mathcal{P}(D)$ and all $\mathcal{Q}_i(D)$ are polynomials in the differential operator, $D$. 

Equivalently, we may simply write,
\begin{equation}
\mathcal{P}(D) x_n(t) = \phi(t) 
\label{eq:gaussccLODAEsEnd}
\end{equation}

where $\mathcal{P}(D)$ is a polynomial of differential operators of equal or higher degree than the individual constitutive equations and $\phi(t)$ is an effective forcing term. This equation is clearly in the form of a higher order governing uni-ccLODE, thereby achieving objective II for sys-ccLODAEs.

Note that if we only need to determine $\mathcal{P}(D)$ - which is the case of the homogeneous sys-ccLODAEs, then we may utilize the constant coefficient nature and alternatively derive it using $\mathcal{P}(D) =$ det$M(D)$. The relationship to the the method described above can be clearly inferred from the fact that the determinant is the multiplication of pivots which are in turn also involved in Gaussian elimination and multiplication of denominators of the $n^{th}$ rows of the elimination matrices.

Because the proposed method for objective II builds on fundamental concepts of linear algebra, existing algorithms may be utilized for symbolic code implementations.

\subsection{Specific Forms of Governing Equation for Small Systems}

The methods for sys-ccLODAEs are suitable for systems of any size. However, many systems, including classically studied systems, have two or three dependent variables and corresponding constitutive equations. We therefore include the specific form for the corresponding higher order governing equations (Equation \ref{eq:gaussccLODAEsEnd}) here. 

\subsubsection{Systems with Two Dependent Variables}

For a 2x2 sys-ccLODAE, the result of Gaussian elimination (Equation \ref{eq:gaussccLODAEs}) is,

\begin{equation}
\bigg( P_{22}(D) - \frac{P_{21}(D)}{P_{11}(D)} P_{12}(D) \bigg) x_2(t) =   \frac{-P_{21}(D)}{P_{11}(D)} f_1(t) + f_2(t)
\end{equation}

Note that for the case of 2x2 sys-ccLODAE this gives an integro-differential equation by use of \ref{eq:partialFractionExpansion} or \ref{eq:factor}. Multiplication of both sides by $P_{11}(D)$, the denominator of last ($2^{nd}$) row of the elimination matrix involved in Gaussian elimination to generate $h(\tilde{P})$, yields the higher order governing equation for $x_2$. Equation \ref{eq:gaussccLODAEsEnd} for 2x2 system is,

\begin{equation}
\underbrace{\bigg( P_{11}(D) P_{22}(D) - P_{21}(D) P_{12}(D) \bigg)}_{\mathcal{P}(D)} x_2(t) =   \underbrace{-P_{21}(D) f_1(t) + P_{11}(D) f_2(t)}_{\phi(t)}
\label{eq:gaussccLODAEsEnd2}
\end{equation}

The proposed method is in fact a systemization of the elimination methods done by hand to derive the governing equation for 2x2 systems. 

\subsubsection{Systems with Three Dependent Variables}

For the three variable system, we get for Equation \ref{eq:gaussccLODAEs}.

\begin{equation}
    h(\Ptilde) = P_{33} - \frac{P_{31}}{P_{11}} P_{13} - \frac{P_{32} - \frac{P_{31}}{P_{11}} P_{12}}{P_{22} - \frac{P_{21}}{P_{11}} P_{12}} \big( P_{23} - \frac{P_{21}
}{P_{11}} P_{13} \big)
\end{equation}

with $P_{ij}$ being polynomials $P_{ij}(D)$. We have dropped explicitly writing $D$ for compactness.

% \begin{equation}
%     h(P) = P_{33}(D) - \frac{P_{31}(D)}{P_{11}(D)} P_{13}(D) - \frac{P_{32}(D) - \frac{P_{31}(D)}{P_{11}(D)} P_{12}(D)}{P_{22}(D) - \frac{P_{21}(D)}{P_{11}(D)} P_{12}(D)} \big( P_{23}(D) - \frac{P_{21}(D)
% }{P_{11}(D)} P_{13}(D) \big)
% \end{equation}

and,
\begin{equation}
    \sum_{i=1}^n g_i(\Ptilde) f_i(t) = \Bigg( -\frac{P_{31}}{P_{11}} -\frac{P_{32} - \frac{P_{31}}{P_{11}} P_{12}}{P_{22} - \frac{P_{21}}{P_{11}} P_{12}} \bigg( \frac{-P_{21}}{P_{11}} \bigg) \Bigg)  f_1 -\frac{P_{32} - \frac{P_{31}}{P_{11}} P_{12}}{P_{22} - \frac{P_{21}}{P_{11}} P_{12}} f_2 + f_3
\end{equation}

where $f_i = f_i(t)$ and we have again dropped explicitly writing the dependence on $t$ for compactness.

By multiplication by denominators of the third rows of the elimination matrices used in Gaussian elimination, $P_{11} \big( P_{22} - \frac{P_{21}}{P_{11}} P_{12} \big)$, we get for Equation \ref{eq:gaussccLODAEsEnd}, the (potentially higher order) governing uni-ccLODE,

\begin{equation}
    \mathcal{P}(D) = (P_{11} P_{22} - P_{12} P_{21}) P_{33 } + P_{12} P_{23} P_{31} + P_{13} (P_{21} P_{32} - P_{22} P_{31} ) - P_{11} P_{23} P_{32}
    \label{eq:gaussccLODAEsEnd3}
\end{equation}

which may also have been obtained directly using the determinant approach.

And,
\begin{equation}
\phi(t) = \big( P_{21} P_{32} - P_{22} P_{31} \big)  f_1(t) - (P_{11} P_{32} -  P_{31}  P_{12}) f_2(t) + \big( P_{11} P_{22} - P_{21} P_{12}  \big) f_3(t)
\end{equation}

% where we have dropped the dependence of the polynomials on the operator $D$ for compactness.

\clearpageSamiya

\subsection{Example: Continuous Coupled Transmission Lines in Frequency Domain}
\label{s:2TL}
 
In this section, we illustrate the method for deriving the governing equation from sys-ccLODAEs using an example motivated by classic electronic traveling-wave amplifiers from microwave theory and an early phenomenological model of cochlear amplification \cite{Hubbard,HubbardRef11}. This system consists of two coupled continuous transmission lines and we use the constitutive equations in the frequency domain as an example. Deriving the governing equation is useful in determining the class of possible behaviors, characterizing the form, and parametric regions of solutions of interest for such a system e.g. regarding power transfer between the lines and determining which parameter choices provide desired amplification patterns.  
% A subset of these purposes {\color{red}(XXXX all ??? Ask Baba)} can be accomplished using the associated characteristic equation for the natural modes.
% which may be useful in determining parameters for design purposes that provide desired amplification patterns.

% NOTE: TO DO: XXXX TO THINK ABOUT: can't you just determine the characteristic equation (if there are no boundary conditions) by setting D to ik in the constitutive equation system, and then you'd just end up with an algebraic equation that you can apply regular methods to to arrive at the characteristic equation... What is advantage: boundary conditions...? TO THINK ABOUT

%In this particular example, the constants $Z(\w), Y(\w)$ are impedances and admittances per unit length NOT POSSIBLE for this to be the case for all of them!

We choose the circuit configuration drawn in \cite{Hubbard} and express constitutive equations relating the dependent variables accordingly. In this system, the top line has variables $I_t(x,\w), V_t(x,\w)$, and the bottom line has variables, $I_b(x,\w), V_b(x,\w)$ \footnote{The identity of the specific parameters and variables are not important - as opposed to the form of these relationships to illustrate the systematic method of deriving the governing equation. In the case of the chosen configuration, $Z_{st}(\w) = j\w L_{ot}, Z_{sb}(\w) = j\w L_{ob}$, $Y_{pt}(\w) = j\w L_t + R_t + \frac{1}{j\w C_t}, Z_{pb}(\w) = j\w C_b + \frac{1}{R_b}$. The effective sources from both lines are due to the voltage across the capacitance.}. We then express these in the form of Equation \ref{eq:ccLODAEsForm} as,

\begin{equation}
\begin{bmatrix} Z_{st}(\w) & D &  & \\ Z_{pt}(\w) D & 1 &  & -\alpha \\  &  & Z_{sb}(\w) & D \\ \beta(\w) D & & D & Y_{bp}(\w) \end{bmatrix}
\begin{bmatrix} I_t(x,\w) \\ V_t(x,\w) \\ I_b(x,\w) \\ V_b(x,\w) \end{bmatrix}
=
0
\end{equation}

where the coupling between the two transmission lines is due to $\alpha$ and $\beta(\w)$, and $D \mapsto \frac{d}{dx}$ is the spatial derivative operator. Notice that setting the coupling parameters to zero leads to separation into two block matrices and, in that case, carrying out the methods of Section \ref{s:governingMethodPart} leads to two independent wave equations (one for the dependent variables of the top line and another for the dependent variables of the bottom line). 

Applying the method described in Section \ref{s:governingMethodPart} with $V_b(\w)$ as the dependent variable of interest \footnote{In this case, the governing equation is the same for all dependent variables because the system is homogeneous and constant coefficient. If the coupling parameters are set to zero, the lowest-order governing equation is different for each of the lines.}, gives after Gaussian elimination (Equation \ref{eq:gaussccLODAEs}),
\begin{equation}
\bigg( Y_{bp} - \frac{\alpha \beta(\w) D^2}{Z_{st} - Z_{pt} D^2} - \frac{D^2}{Z_{sb}} \bigg) V_b(x,\w) = 0
\end{equation}

Multiplication by the relevant denominators of the elimination matrices $Z_{sb}(\w) \bigg( Z_{st}(\w) - Z_{pt}(\w) D^2  \bigg)$ \footnote{We also divide by $Z_{pt}(\w)$ to put the governing equation in the standard ODE form in which the coefficient of the highest order differential operator is unity.}, yields the higher order governing equation (Equation \ref{eq:gaussccLODAEsEnd}), 

\begin{equation}
\Bigg(D^4 - \bigg( \frac{Z_{st} (\w)}{Z_{pt} (\w)} + \alpha \beta(\w) \frac{Z_{sb} (\w)}{Z_{pt} (\w)} + Z_{sb}(\w) Y_{bp}(\w) \bigg) D^2  + \frac{Z_{st} (\w) Z_{sb} (\w) Y_{bp} (\w)}{Z_{pt} (\w)} \Bigg) V_b(x,\w) = 0
\end{equation}

which is conducive to studying the system. Note that we may have arrived at this equation directly from the sys-ccLODAEs using $\mathcal{P}(D) = $det$M(D)$ as described in Section \ref{s:governingMethodPart} because the system is homogeneous, but sought to exemplify the method for the more general case. 

To return to the block diagonal case: in the case of  $\alpha = \beta = 0$, the $\mathcal{P}(D)$ of the above equation may be factored into two second order differential operator polynomials - or two wave equations with wavenumbers $\sqrt{\frac{Z_{st} (\w)}{Z_{pt} (\w)}}, \sqrt{Z_{sb} (\w) Y_{sb} (\w)}$. Each of which provides two solutions. The boundary conditions impose that one set of solutions describes behavior on one transmission line and the other set describes behavior on the other transmission line.

\clearpageSamiya

\subsection{Example: Continuous Transmission Line Including Algebraic Equation in Frequency Domain}

% In this section, we apply our method for deriving the governing equation to a system of ccLODAEs that includes an algebraic equation to illustrate that our method is not restricted to systems of ccLODEs.

In this section, we illustrate that our method for deriving the governing equation applies not only to systems of sys-ccLODEs but also to sys-ccLODAEs which also occur naturally in systems and phenomenon - though they are less-often the focus of development of even numerical methods. The sys-ccLODAE is the constant coefficient version of a classical model of auditory filters - e.g. \cite{zweig76}, that may be appropriate for design purposes (such as machine hearing) rather than a physiological model of the cochlea. We define $D \mapsto \frac{d}{dx}$ and choose the variable of interest for this sys-ccLODAE to be $V_{bm}(x,\w)$. Expressed as Equation \ref{eq:ccLODAEsForm}, the system is,

\begin{equation}
    \begin{bmatrix}
    D & Z(\w) & \\
    & D & b \\
    1 & & b Z_{fbal}(\w)
    \end{bmatrix}
    \begin{bmatrix}
    P(x,\w) \\ U(x,\w) \\ V_{bm}(x,\w)
    \end{bmatrix}
    = 
    \begin{bmatrix}
    0\\0\\0
    \end{bmatrix}
\end{equation}

From Equation \ref{eq:gaussccLODAEsEnd3} (and division by $b Z_{fbal}$), we obtain the governing equation,
\begin{equation}
    \bigg(D^2 + \frac{Z}{Z_{fbal}} \bigg) V_{bm}(x,\w) = 0
\end{equation}

\clearpageSamiya

\section{Methods for Solving sys-ccLODAEs} 
\label{s:solutionMethod}
% $https://en.wikipedia.org/wiki/Symbolic_integration$
In this section, we develop methods for objective III - deriving analytic expressions for the full solutions for a system of multivariate constitutive equations that form a sys-ccLODAE. The problem is solved in the native time or space domain of the problem. In \ref{s:differentialOperator}, we review differential operator methods for analytically solving uni-ccLODEs and make some modifications. Those methods are utilized as part of our approach for solving sys-ccLODAEs.

% (as opposed to first requiring deriving the governing uni-ccLODEs for each dependent variable which is increasingly difficult for larger systems). The solutions are obtained in the native domains - as opposed to, for instance, using Laplace transform with the singular state-space formulation (Equation \ref{eq: stateSpaceDAE}) to solve for the state variables and which involves imposed limitation to IVPs.

\subsection{Method}

For objective III, we again use concepts from linear algebra. Recall that to solve a sys-LAE as in Equation \ref{eq:sysLAE}, we use inversion of the nonsingular matrix $A$, which results in solutions for all $x_i$,
\begin{equation}
    x_i = \sum_{j=1}^{n} h_{j;i}(l) b_j
    \label{eq:inversionResultLinAlg}
\end{equation}

The Gauss-Jordan method for inversion takes an algebraic system of equations, applies row operations to both sides, up to arriving at pivot columns of $A$. For the case of the full-rank matrix $A$, after row-normalization, this leads to the identity matrix, which is then used to determine the solution of each of the variables in the column vector $x$. Inversion, and particularly the Gauss-Jordan method is relevant for determining the full solutions for each of the dependent variables in a system of sys-LDAEs. The resultant expressions are analytic (but not necessarily closed-form).

The key to developing our methods for solving sys-ccLODAEs is to formulate the sys-ccLODAEs as in Equation \ref{eq:ccLODAEsForm}. If the only interest is in the complementary solutions (which are the same for all dependent variables for sys-ccLODAEs) or the system is homogeneous, then the solution can be obtained by determining the roots of the characteristic equation,
\begin{equation}
    \textrm{det}M(a) = 0
\end{equation}

If we are interested in the full solution (including complementary solutions and a particular) to a nonhomogeneous sys-ccLODAE, then we invert the matrix of operators. Due to the constant nature of the coefficients, we may define the inverse of the matrix of operators $M(D)$ by extrapolating standard matrix inversion methods.

If the inverse is obtained using Gauss-Jordan elimination, then these operations result in an equation parallel in form to Equation \ref{eq:inversionResultLinAlg},

\begin{equation}
    x_i(t) = \sum_{j=1}^n h_{j;i}(D) f_j(t)
    \label{eq:formSolBefore}
\end{equation}

where $h_{j;i}(D)$ has a form due to the row operations for inversion described by Gauss-Jordan elimination. We then multiply each $h_{j;i}$ by fractions of equal numerator and denominator (unity) repeatedly to put $x_i(t)$ in the form,

\begin{equation}
    x_i(t) = \sum_{j=1}^{n} \frac{P_{j;i}(D)}{Q_{j;i}(D)} f_j(t) = \sum_{j=1}^{n}\frac{1}{Q_{j;i}(D)}\phi_{j;i}(t)
    \label{eq:formSol}
\end{equation}

which can then be solved for each of the dependent variables using differential operator methods for uni-ccLODEs described in \ref{s:differentialOperator} (Equation \ref{eq:partialFractionExpansion} or \ref{eq:factor}) with Definition \ref{eq:modifiedDefinitionIntegral} for the full solution or Definition \ref{eq:classicalDefinitionIntegral} for the particular solutions only.

The methods are appropriate for sys-ccLODAEs of any size, and are particularly simple and intuitive due to their relation to standard matrix inversion (which is possible due to the constant coefficient nature of the equations). Because the methods we develop in this paper build on fundamental concepts of linear algebra, existing algorithms may be utilized for symbolic code implementations. However, we note that methods for objective II only requires elimination and are therefore more computationally efficient than those for objective III which involves matrix inversion.

The methods for objective III provide analytic solutions which are useful in studying the system in ways that cannot be achieved if we are limited to numerical simulation of systems of differential equations which suffer from issues such as accuracy, stability, specificity to parameters, boundary and initial conditions, and discretization schemes, and limitations to types of equations. Note that the solution method does not require transforms to move to a different domain but rather solves the problem in its native domains. The analytic method is not preferential to any particular type of problem (e.g. IVPs or BVPs).

\clearpageSamiya
% \subsection{Specific Form of Solution for Small Systems}

%  Again, as many classical systems are small, we provide the solution for 2x2 and 3x3 systems. For system with two dependent variables, the solution in the form of Equation \ref{eq:formSol} (easily seen in relation to the inverse of algebraic matrices) is,

% \begin{equation}
%     x(t) = \frac{1}{P_{11}(D)P_{22}(D)- P_{12}(D)P_{21}(D)}
%     \begin{bmatrix}
%     P_{22}(D) f_1(t) & -P_{12}(D) f_2(t) \\ -P_{21}(D) f_1(t) & P_{11}(D)f_2(t)
%     \end{bmatrix} 
%     \label{eq:solTwoByTwo}
% \end{equation}

% and for a 3x3 system,
% \begin{equation}
%     x(t) = 
% \end{equation}

% where the dependence on $D$ and $t$ is not written explicitly for compactness. 

% XXXXX TO DO: XXXX TO THINK ABOUT XXXX

\clearpageSamiya
\subsection{Example: Forced Vibration Systems}

To exemplify our method for solving sys-ccLODAEs, consider a typical mass-spring system with two masses $m_1, m_2$, connected via a spring with constant $k_1$. The second mass is connected to a fixed wall with a spring with constant $k_2$. An external upward force $f_1(t)$ is applied to the $m_1$, and the upward displacements of the two bodies are $x_1(t), x_2(t)$.

The system of constitutive equations is,
\begin{equation}
    \begin{cases}
    m_1 \ddot{x}_1(t) + k_1 x_1(t) - k_1 x_2(t) = f_1(t) \\
    m_2 \ddot{x}_2(t) - k_1 x_1(t) + (k_1 + k_2) x_2(t) = 0
    \end{cases}
\end{equation}

Our proposed method for solving sys-ccLODAEs is applicable to vibration systems of any configuration and admits forcing terms. It is not restricted to pure mass-spring systems but can include dashpots. To solve for $x_1(t), x_2(t)$ using our method, we first formulate the system as in Equation \ref{eq:ccLODAEsFormDetails}. With $D \mapsto \frac{d}{dt}$, we express this system as,

\begin{equation}
   \begin{bmatrix}
    m_1 D^2 + k_1 & -k_1 \\
    -k_1 & m_2 D^2 + k_1 + k_2
    \end{bmatrix} 
    \begin{bmatrix}
    x_1(t) \\ x_2(t)
    \end{bmatrix}
    = 
    \begin{bmatrix}
    f_1(t) \\ 0
    \end{bmatrix}
\end{equation}

Inverting for $x(t)$ as described in Section \ref{s:solutionMethod}, we obtain,
\begin{equation}
    x(t) = \frac{1}{(m_1 D^2 + k_1)(m_2 D^2 + k_1 + k_2) - k_1^2}
    \begin{bmatrix}
    m_2 D^2 + k_1 + k_2 & k_1\\ k_1 & m_1 D^2 + k_1
    \end{bmatrix} 
    \begin{bmatrix}
    f_1(t) \\ 0
    \end{bmatrix}
\end{equation}

We can solve this for each dependent variable by using the methods described in more detail in \ref{s:differentialOperator}. Factoring the polynomial in the denominator as
\begin{equation}
    \frac{1}{(m_1 D^2 + k_1)(m_2 D^2 + k_1 + k_2) - k_1^2} = 
    \frac{1}{m_1 m_2} \frac{1}{(D-a)(D+a)(D-b)(D+b)} 
\end{equation}

with $a,b = \frac{i}{\sqrt{2}} \sqrt{(\frac{k_1}{m_2} + \frac{k_1}{m_1} + \frac{k_2}{m_2}) \pm \sqrt{ \frac{k_1^2}{m_1^2} + \frac{k_1^2}{m_2^1} + \frac{k_2^2}{m_2^2} + \frac{2 k_1 k_2}{m_2^2} + \frac{2 k_1^2}{m_1 m_2} - \frac{2k_1 k_2}{m_1 m_2}}}$, we then obtain analytic expressions for the full solution for each of the dependent variables $x_1(t), x_2(t)$ \footnote{Note that this does not provide a relationship between the constants of the complementary components of the two dependent variables. }. As there are no repeated first order factors of the polynomial in the denominator, we may use either Equation \ref{eq:partialFractionExpansion} or Equation \ref{eq:factor} to obtain these expressions (only Equation \ref{eq:factor} is appropriate in the case of repeated first order factors  - see \ref{s:differentialOperator}). For example, using Equation \ref{eq:partialFractionExpansion} to solve for $x_2(t)$, we obtain,

\begin{align}
    x_2(t) &= \frac{1}{(m_1 D^2 + k_1)(m_2 D^2 + k_1 + k_2) - k_1^2} k_1 f_1(t) \\
    &=\frac{k_1}{2 m_1 m_2}  \frac{1}{a^2-b^2} \Bigg( \frac{1}{a} \expn{at} \defint{}{}{\expn{-at} f_1(t)}{t}  - \frac{1}{a} \expn{-at} \defint{}{}{\expn{at} f_1(t)}{t} \nonumber \\& - \frac{1}{b} \expn{bt} \defint{}{}{\expn{-bt} f_1(t)}{t}  + \frac{1}{b} \expn{-bt} \defint{}{}{\expn{bt} f_1(t)}{t}  \Bigg)
\end{align}

Equivalently, we may express this using Equation \ref{eq:factor}, 
\begin{align}
x_2(t) &= \frac{1}{(m_1 D^2 + k_1)(m_2 D^2 + k_1 + k_2) - k_1^2} k_1 f_1(t) \\
&= \frac{k_1}{m_1 m_2}\expn{at}\defint{}{}{ \expn{-2at}\defint{}{}{\expn{(a+b)t}\defint{}{}{\expn{-2bt} \defint{}{}{\expn{bt} f_1(t)}{t}}{t}}{t}}{t}
\end{align}

For instance, for $f_1(t) = f_1$, we get the full solution, using either expression,
\begin{equation}
    x_2(t) = \frac{k_1 }{m_1 m_2 a^2b^2}  f_1
    + C_1 \expn{at} + C_2 \expn{-at} + C_3 \expn{bt} + C_4 \expn{-bt}
\end{equation}

which clearly includes the complementary solutions and a particular solution. The method yields both particular and complementary parts of the solution together rather than solving for each separately.

\clearpageSamiya

\section{Methods and Limitations for sys-ccLPDAEs}
\label{s:ccLPDAEs}
In this section, we present our methods for objectives I-III (formulation, governing equation derivation, and system solutions) starting from systems of constitutive equations that are not sys-ccLODAEs but rather sys-ccLPDAEs where we have multiple independent variables and the corresponding governing equation is a uni-ccLPDE as opposed to a uni-ccLODE. In order to extend our methods and determine their applicability to, and limitation for, sys-ccLPDAEs, we must first extend differential operator theory to the multi-operator case.

\subsection{Differential Operator Theory for uni-ccLPDEs}

We may define each differential operator $D_{i}$ to be associated with an independent variable $t_i$, such that,
\begin{equation}
    D_i x(t_1, t_2, \dots t_m)  \triangleq \frac{\partial}{\partial t_i}x(t_1, t_2, \dots t_m)
\end{equation}

Consequently, we may generally express a uni-ccLPDE as,

\begin{equation}
\mathcal{P}(D_1, D_2, \dots D_m) x(t_1, t_2, \dots t_m) = \phi (t_1, \dots t_m) 
\label{eq:ccLPDEform}
\end{equation}

with $\mathcal{P}(D_1, \dots D_m)$ being a polynomial in the differential operators.

The multi-operator counterparts to the properties of \ref{s:differentialOperator} apply for any $D_i$.

Accordingly, for instance, 
\begin{equation}
P(D_1, D_2, \dots D_m) \expn{a_1 t_1 + a_2 t_2 + \dots a_m t_m} = P(a_1, a_2, \dots a_m) \expn{a_1 t_1 + a_2 t_2 + \dots a_m t_m} 
\end{equation}

Note that this property can be used to simply find the relationship between the different components of the complementary solution to a uni-ccLPDE by setting $P(a_1, a_2, \dots a_m)=0$. This is the counterpart to the characteristic equation obtained from a uni-ccLODE (see Equation \ref{eq:toCharacteristicEq}). 

Similar to Equation \ref{eq:interchangeable}, the order of integration and differentiation does not matter when applied to smooth functions regardless of the independent variables and whether they are the same or mixed independent variables. This is due to the constant coefficient nature of the system.

\clearpageSamiya

\subsection{Formulating Constitutive sys-ccLPDAEs as Linear System of Equations}

To formulate the sys-ccLPDAEs (objective I), we express the system of equations as in Equation \ref{eq:ccLODAEsForm}, but with $\mathbf{t}$ a vector of $n$ independent variables $t_1, t_2, \cdots t_n$ (e.g. time and spatial coordinates), and $\mathbf{D}$ a vector of differential operators. Accordingly, the uni-operator polynomials $P_{ij}(D)$ of sys-ccLODAEs are replaced by multi-operator polynomials in multiple differential operators $P_{ij}(\mathbf{D}) = P_{ij}(D_{t_1}, D_{t_2}, \cdots D_{t_n})$ for sys-ccLPDAEs. We therefore express the counterpart of Equation \ref{eq:ccLODAEsForm} for sys-ccLPDAEs as,

\begin{equation}
    M(D_1, D_2, \cdots D_{t_n}) x(t_1, t_2, \cdots t_n) = f(t_1, t_2, \cdots t_n)
\end{equation}

or, compactly,
\begin{equation}
    M(\mathbf{D}) x(\mathbf{t}) = f(\mathbf{t})
    \label{eq:ccLPDAEsForm}
\end{equation}

with $M(\mathbf{D})$ the matrix of operators, $x(\mathbf{t})$ the vector of dependent variables, and $f(\mathbf{t})$ the vector of forcing terms.

\clearpageSamiya
\subsection{Derivation of Governing uni-ccLPDE from Constitutive sys-ccLPDAEs}

\subsubsection{Method}

The method for deriving higher order governing equations (objective II) for constitutive sys-ccLPDAEs, transfers directly from sys-ccLODAEs (Section \ref{s:governingMethod}) with no limitations or approximations but with the aforementioned changes to the independent variables and operators.

\subsubsection{Example: Continuous Transmission Lines in Time Domain}

In this section, we exemplify the method for deriving the higher order governing equation from a system of constitutive sys-ccLPDAEs. For our example, we consider a transmission line (which is of interest in many propagation systems \cite{speechTL, BabaDistributedSpeech, cochleaTL}) in the time domain. A commonly encountered transmission line has an inductor with inductance $L_1 \delta x$ for the series component, and a simple resonant harmonic oscillator or RLC in series for the parallel component \footnote{Trivially, setting $R = L = 0$ leads to the simple 1D wave equation as becomes apparent in Equation \ref{eq:waveEq}. Addition of nonzero $R,L$ adds dissipative (or, depending on the sign, growing) and dispersive behavior.}. The constitutive equations arise from Kirchhoff's current law (continuity) and Kirchoff's voltage law (force balance), and the forcing terms are chosen to be zero in this example. The system responds to a set of initial and boundary conditions. The differential equations are defined by taking the limit of differential equations as $\delta x \xrightarrow{} 0$. The set of constitutive equation for this system is,

\begin{equation}
    \begin{cases} 
    L \frac{\partial^2}{\partial t \partial x}  i(x,t) + R \frac{\partial}{\partial x}  i(x,t)  + \frac{1}{C} \defint{-\infty}{t}{\frac{\partial}{\partial x}  i(x,\tau)}{\tau} + v(x,t) = 0 \\
    L_1 \frac{\partial}{\partial t} i(x,t) + \frac{\partial}{\partial x}  v(x,t) = 0
    \end{cases}
\end{equation}

We use $D_t$ and $D_x$ to denote time and space differential operators, $D_t \mapsto \frac{\partial}{\partial t}, D_x \mapsto \frac{\partial}{\partial x}$. We express the constitutive equations as a linear system. This is a system in which the integro-differential equation can be reformulated as a differential equation simply by multiplication of the first row by $D_t$ to use the methods developed in this paper. This yields for objective I (Equation \ref{eq:ccLPDAEsForm}),

% \begin{equation}
%     \begin{bmatrix}
%   ( R + \frac{1}{C} \frac{1}{D_{t}} +  L D_{t} )D_{x} & 1\\   L_{1} D_{t} & D_{x} \end{bmatrix}
% \begin{bmatrix} i(x,t) \\ v(x,t) \end{bmatrix}
% =
% \begin{bmatrix} 0 \\ 0 \end{bmatrix}
% \end{equation}

\begin{equation}
    \begin{bmatrix}
   (L D^2_{t} +  R D_t + \frac{1}{C} )D_{x} & 1\\   L_{1} D_{t} & D_{x} \end{bmatrix}
\begin{bmatrix} i(x,t) \\ v(x,t) \end{bmatrix}
=
\begin{bmatrix} 0 \\ 0 \end{bmatrix}
\end{equation}

% Note that the integral component in the first row is due to the capacitor. 

We then apply our method to derive the governing equation (objective II) for $v(x,t)$. Using Equation \ref{eq:gaussccLODAEsEnd2}, we obtain,

\begin{equation}
\bigg(\big(LD^2_t + R D_t + \frac{1}{C} \big) D^2_{x}  - L_1 D^2_t \bigg) v(x,t) = 0
\label{eq:waveEq}
\end{equation}

which can then be characterized according to its parameter values which may be chosen or controlled for design problems.

As mentioned previously, derivation of the governing equation allows us to characterize the system and determine parametric conditions for various behaviors without the need for solving the constitutive equations - e.g. conditions under which wave propagation is supported and form of exchange of energy between the coupled two transmission line system of Section \ref{s:2TL}. 
%{\color{red} XXXX why can't they just do this with characteristic equations after moving to the w-k domain. Ask Baba?!}.

\clearpageSamiya

\subsection{Solving for Dependent Variables: Methods, Approximations, and Limitations}

In this section, we discuss methods and limitations for objective III - obtaining exact or approximate analytic solution, for sys-ccLPDAEs. We first formulate the sys-ccLPDAEs as in Equation \ref{eq:ccLPDAEsForm} and then carry out the methods of Section \ref{s:solutionMethod} to arrive at the multi-operator counterpart to Equation \ref{eq:formSol}. This results in, for each dependent variable $x_i(\mathbf{t})$,

\begin{equation}
    x_i(\mathbf{t}) = \sum_{j=1}^{n}\frac{1}{Q_{j;i}(\mathbf{D})}\phi_{j;i}(\mathbf{t})
    \label{eq:formSolccLPDAEs}
\end{equation}

This equation cannot generally be solved using methods of differential operator theory. The issue is not due to  dealing with a system, but is rather inherent in partial differential equations. To our knowledge, there are no general analytic solutions for nonhomogeneous uni-ccLPDEs (Equation \ref{eq:ccLPDEform}) \footnote{This can be seen by the fact that unlike $\frac{1}{P(D)}$ of \ref{s:differentialOperator},  $\frac{1}{P(D_1, D_2, \cdots, D_n)}$ is not defined as an operator. For instance, we cannot solve an equation $(D_x + D_t) u(x,t) = f(x,t)$, by simply expressing it as, $u(x,t) = \frac{1}{D_x + D_t}f(x,t)$}. Certain classical uni-ccLPDEs (predominantly homogeneous and without mixed derivatives) can be solved using the method of separation of variables.    
% {\color{red}  How does your method relate to separable solutions?}

We may solve a uni-ccLPDE if we can, exactly or approximately, factor the multi-operator polynomial into uni-operator polynomials (which may require transformations) -  i.e. we can express  $\frac{1}{P(D_1, D_2, \cdots, D_n)} \approx \frac{1}{P_1(\tilde{D}_1)} \frac{1}{P_2(\tilde{D}_2)} \cdots \frac{1}{P_m(\tilde{D}_m)}$. This allows us to then apply equations \ref{eq:partialFractionExpansion} or \ref{eq:factor} to arrive at the analytic solution for such a uni-ccLPDE.

Consequently, we may express the analytic solutions for sys-ccLPDAEs for cases of naturally weak coupling in the native independent variables as,
\begin{equation}
    x_i(\mathbf{t}) \approx \sum_{j=1}^{n}\prod\limits_{k=1}^{m}\frac{1}{Q_{j;i;k}(D_k)}\phi_{j;i}(\mathbf{t})
\end{equation}

which can then be solved analytically using equations \ref{eq:partialFractionExpansion} or \ref{eq:factor}. This requires factoring the multi-operator polynomials $Q_{j;i}(\mathbf{D})$ of Equation \ref{eq:formSolccLPDAEs} such that the absolutely irreducible factors are uni-operator polynomials (i.e. each factor only contains $D_k$ for a single $k$).

Only a subset of sys-ccLPDAEs naturally have weak coupling in their independent variables. Therefore, in order to factor $Q_{j;i}(\mathbf{D})$ into uni-operator polynomials, we may additionally require manipulating the equation using transformations of the independent and/or dependent variables.

\begin{equation}
    \tilde{x}_i(\mathbf{\tilde{t}}) \approx \sum_{j=1}^{N}\prod\limits_{k=1}^{M}\frac{1}{\tilde{Q}_{j;i;k}(\tilde{D}_k)}\tilde{\phi}_{j;i}(\mathbf{\tilde{t}})
\end{equation}

%Note that the factorization can be applied to a modified equation with known transformations mapping . 

%The transformation map is chosen to approximately decouple the components of the corresponding homogeneous solution. 

%  In such an approximate space, the dependence of the dependent variable on the various independent variables is essentially decoupled

% Condition is the for integration (denominators) multiplication of polynomials of different operators. This condition is equivalent to saying that the components of the corresponding homogeneous solutions (in some transformed independent variables) are weakly coupled so that their coupling may be neglected.

% XXXX TO DO: XXXXX check if differential operator methods already exist for ccLPDEs XXXXX

\clearpageSamiya

\section{On sys-vcLODAEs and sys-vcLPDAEs}
\label{s:vc}

In this section, we present a formulation for sys-vcLODAEs which is also extendable to sys-vcLPDAEs (objective I). We then consider special cases of sys-vcLODAEs for which we can derive the governing equation (objective II) using the same methods as sys-ccLODAEs. This also applies to the corresponding special cases of sys-vcLPDAEs. In this paper, we do not discuss methods for determining solutions to sys-vcLODAEs and sys-vcLPDAEs, and it is \textit{generally} not possible to provide exact analytic solutions for these categories.

\subsection{Formulating the System of Equations}
We may formulate the sys-vcLODAEs as follows,

\begin{equation}
    \mathcal{M}(D_t, t) x(t) = f(t)
    \label{eq:formulationvcLODAEs}
\end{equation}

with $D_t \mapsto \frac{d}{dt}$. The elements of the matrix $\mathcal{M}$ are $H_{ij}(D_t, t)$, which take the form,

\begin{equation}
     H(D_t, t) = \sum \limits_{k=0}^{n} a_k(t) D_t^k
\end{equation}

Note that sys-ccLODAEs are special cases of sys-vcLODAEs when the matrix elements $H_{ij}(D_t, t) = P_{ij}(D_t)$.

As a particularly simple example, consider a continuous transmission line with incremental sources and in which the series components are constant along $x$ but the parallel components vary with $x$. The system of equations (sys-vcLODAEs) is formulated as,

\begin{equation}
    \begin{bmatrix}
    Z(\w) & D \\ D & Y(x,\w)
    \end{bmatrix}
    \begin{bmatrix}
    I(x,\w) \\ V(x,\w)
    \end{bmatrix} = 
    \begin{bmatrix}
    f_{1}(x,\w) \\ f_{2}(x,\w)
    \end{bmatrix}
    \label{eq:vcExample}
\end{equation}

\subsection{Deriving the Governing Equation for Special Cases of sys-vcLODAEs and sys-vcLPDAEs}
\label{s:vcGoverning}

\subsubsection{Method}

In this section, we discuss methods for deriving the governing equation (objective II) for special cases of sys-vcLODAEs and sys-vcLPDAEs. We cannot generally extrapolate the methods we developed for deriving the governing equation for sys-ccLODAEs to the case of sys-vcLODAEs. Two steps of the methods for sys-ccLODAEs account for this: 
\begin{enumerate}
    \item The order of operations of Gaussian elimination to arrive at Equation \ref{eq:gaussccLODAEs} does not matter in the case of constant coefficients, but does matter for variable coefficients.
    \item Multiplication of Equation \ref{eq:gaussccLODAEs} by denominators results in the governing uni-ccLODE (Equation \ref{eq:gaussccLODAEsEnd}) for sys-ccLODAEs because multiplication and division are interchangable for the sys-ccLODAEs. However this is not the case for general sys-vcLODAEs.
\end{enumerate}

There is a subset of sys-vcLODAEs (Equation \ref{eq:formulationvcLODAEs}) for which these two issues do not arise. For these cases, we simply use the methods we constructed for sys-ccLODAEs (Section \ref{s:governingMethod}) for deriving the corresponding governing equation. The criteria is that the sys-vcLODAEs are of the form,

\begin{equation}
\begin{bmatrix}
P_{1,1}(D_t) & \cdots & P_{1,n-1}(D_t) & H_{1,n}(D_t,t)\\
\vdots &  \ddots & \vdots & \vdots \\
P_{n,1}(D_t) & \cdots & P_{n,n-1}(D_t) & H_{n,n}(D_t,t)
\end{bmatrix} 
\begin{bmatrix} x_1(t) \\ \vdots \\ x_n(t) \end{bmatrix}
= 
\begin{bmatrix} f_1(t) \\ \vdots \\ f_n(t) \end{bmatrix} 
    \label{eq:lastColumn}
\end{equation}

In other words, only the final column - operating on $x_n(t)$, contains the variable coefficients. In other words, all variable coefficients in the system of equations must correspond to the same dependent variable. For such systems, issues (1) and (2) do not arise because the corresponding elimination matrices used for Gaussian elimination contain only constant coefficient terms \footnote{Equivalently, the $L$ in $LU$ matrix decomposition contains only constant coefficient terms} (the elements in the final column of $\mathcal{M}$ do not appear as part of the elimination matrices). 

As for sys-vcLPDAEs, the governing equation is derived in a similar manner and under the same condition as sys-vcLODAEs. Future directions may include methods for deriving the governing equation for more general cases of sys-vcLODAEs and sys-vcLPDAEs in which more than one dependent variable has variable coefficients. Future directions also include developing methods for deriving approximate expressions for the solution for each dependent variable in a sys-vcLODAE.

\subsubsection{Example: Continuous Transmission Line with Varying Parallel Admittance}

To give an example, we use the system in Equation \ref{eq:vcExample} and choose $V(x, \w)$ to be the dependent variable of the governing equation. Because it satisfies the condition in Equation \ref{eq:lastColumn}, we may express the governing equation using Equation \ref{eq:gaussccLODAEsEnd2} as,

\begin{equation}
    \bigg( D_{x}^2 - Z(\w) Y(x,\w) \bigg) V(x,\w) = D_{x} f_1(x,\w) - Z(\w) f_2(x,\w)
\end{equation}

Note that the condition for extrapolating sys-ccLODAE methods for the governing equation is not met for the same system if the dependent variable of interest is $I(x,\w)$. Consequently, if given flexibility regarding the choice of dependent variable with regards to governing equations, it is advantageous to choose dependent variables for which the condition (Equation \ref{eq:lastColumn}) holds.

% \subsection{Further Directions}
% Future directions may include developing approximate methods for deriving the governing equation for the more general case of vcLODAEs and deriving analytic expressions for the full solutions for the system of equations. 

\clearpageSamiya

\section{Conclusion}

\subsection{Summary}

Analytic solutions to the system and analytic characterization through governing equations provide insights into the behaviors of DAEs as well as the parametric regions of operation for each potential behavior. This is relevant for interpretable modeling, analytic characterization and estimation, and engineering design in which the objective is to tune parameter values to achieve specific behavior. Such insights cannot directly be obtained using numerical simulations. In this paper, we introduced methods for deriving analytic solutions from DAEs as well as methods for deriving governing equations for analytic characterization.

We first introduced a direct formulation for representing systems of linear differential-algebraic scalar constitutive equations that are homogeneous or nonhomogeneous with constant or nonconstant coefficients. Using this formulation for the system of constitutive equations (DAEs), we (a) developed simple methods for deriving analytic expressions for the full solution (complementary and particular) for all dependent variables of systems that consist of constant coefficient ordinary-DAEs and special cases of partial-DAEs. We also (b) developed simple methods for deriving the governing equation for a chosen dependent variable for the constant coefficient ordinary-DAEs and partial-DAEs as well as special cases of variable coefficient DAEs. 

The proposed formulation to represent systems of DAEs uses differential operators. Unlike state space formulations, our formulation follows very directly from the system of constitutive equations, directly handles algebraic equations, handles higher order derivatives without requiring defining new variables, and is appropriate for systems containing PDEs. To develop the methods for deriving the governing equation and the analytic solutions from the system of constitutive equations, we utilized the concepts of Gaussian elimination and inversion from linear algebra and (somewhat modified) differential operator theory. 

The methods are simple, build on fundamental linear algebra and hence existing algorithms may be used, can be automated with symbolic coding environments thereby allowing for dealing with systems of any size while retaining analytic nature, and are intended for use with a variety of engineering, physical, and financial systems. We exemplified our methods using dynamic and propagation systems such as vibrating masses and coupled continuous transmission lines.

\subsection{Future Directions}

Future directions include addressing limitations of the methods presented here. This includes developing methods for deriving the exact or approximate governing equation for the general case of sys-vcLDAEs as opposed to only those in which all variable coefficients correspond to a single dependent variable. This also includes methods for deriving approximate analytic expressions for the full solutions for dependent variables in sys-vcLODAEs. Another limitation is regarding solving the system of sys-ccLPDAEs, and future work includes developing approaches for determining transformations and approximate uni-operator factorization methods in order to derive analytic solutions for the  sys-ccLPDAE. Methods for exact and approximate factorization of multivariate polynomials exist in the literature but to our knowledge these do not generally result in univariate factors. Instead, they represent or approximate a multivariate polynomial as a multiplication of multivariate (rather than univariate) polynomial factors \cite{methodsfactoringpolynomial1, methodsfactoringpolynomial2, methodsfactoringpolynomial3, methodsfactoringpolynomial4}. Hence, future directions include developing systematic methods for constructing the transformations of $\mathbf{t}, x_i, \phi_{j;i}$ into $\mathbf{\tilde{t}}, \tilde{x}_i, \tilde{\phi}_{j;i}$ as well as constructing approximate factorization methods. Doing so would enable us to obtain approximate analytic solutions for sys-ccLPDAEs other than those that are weakly coupled in their native independent variables and are for which the factorization is apparent.

Future directions also include extending the framework to develop methods for deriving the governing equation and analytic expressions for full solutions for related classes of problems - specifically linear systems with other operators which share certain properties with differential operators (e.g. linear difference-algebraic equations). % \cite{PseudoLinearAlgebra} 
As our formulation is in a simple matrix form $Mx = f$, we may, in the future, be able to extend linear algebra concepts such as factorization and pseudoinverse to $M$ of operator polynomials. Finally, the formulation may also be used for developing numerical solution methods for boundary value problems.

\appendix

\section{Classical Formulation and Solution Methods}
\label{s:literatureSurvey}

In this appendix, we discuss the main analytic methods for solving some classes of sys-LDAEs. As mentioned previously, we are only concerned with smooth linear systems. 

The class of sys-LODEs are the most studied. They are represented using the state-space formulation, which is of the form,

\begin{equation}
    \dot{x}(t) = A(t)x(t) + B(t)u(t)
    \label{eq:stateSpace1}
\end{equation}

Constitutive equations are not necessarily native to this form. Instead, state variables are defined so that constitutive equations of a system can be formulated as in Equation \ref{eq:stateSpace1}. Note that the formulation is a set of first order ODEs where the derivative is in a single dependent variable. Second order undamped mass-spring systems $M\ddot{x}(t) + Kx(t) = 0$ effectively fall into this category. Current methods for solving or analyzing these sys-LODEs are built on the state-space formulation.

There are several numerical algorithms for solving sys-LODE initial value problems (IVPs) of this formulation. The constant coefficient problems, sys-ccLODE, in the state-space form have simple analytic solutions. In the time-domain, the solution is achieved through matrix exponentials and integrating factors $x(t) =  \expn{At} x(0) + \defint{0}{t}{\expn{A(t-\tau)} B u(\tau)}{\tau}$.
% https://en.wikibooks.org/wiki/Control_Systems/Linear_System_Solutions
Alternatively, the analytic solution can be achieved by Laplace transform $sX(s) -x(0) = AX(s) + BU(s)$ which is then solved as an eigenvalue problem or converted into transfer functions for a single variable and solved accordingly \footnote{Conversion to a transfer function is in some ways similar to deriving a governing equation from the system of constitutive equations. This is because they are in terms of a single dependent variable, and it is possible to convert some governing equations into transfer functions by taking the Laplace transform after deriving the governing equation from the system of constitutive equations}. Note that these are solved in a domain than the independent variables. Laplace transform methods are used for LTI systems and cannot be used for systems with variable coefficients. Beyond analytic solutions for sys-ccLODEs, there is literature on the analysis of such systems from a controls perspective (controllability, stability, solvability, observability). The study of control properties also builds on the state space formulation and is useful for designing observers and controllers - which are \textit{not} the objectives of this paper. 

Variable coefficient sys-LODE IVPs that can be expressed in the form $\dot{x}(t) = A(t) x(t) + B(t) u(t); x(0) = x_o$ have a solution $x(t) = \phi(t,0) x(0) + \defint{0}{t}{\phi(t,\tau)B(\tau) u(\tau)}{\tau}$. However, $\phi(t,\tau)$ is a problem-specific transition matrix that is not always easy to determine. For the case of sys-ccLODEs in the previous paragraph, $\phi(t,\tau) = \expn{A(t-\tau)}$. Special cases of sys-vcLODEs admit closed-form expressions for $\phi(t,\tau)$, but this is generally not the case. The solution to certain homogeneous sys-vcLODE IVPs can be expressed using Magnus expansion which is useful for analytic study when the series converges \cite{Magnus}. The higher terms in this series are increasingly complex, and hence it is desirable to approximate by truncation.

Sys-LODAEs are also formulated using a more general state-space representation,

\begin{equation}
E(t)\dot{x}(t) = A(t)x(t) + B(t)u(t)
 \label{eq:stateSpaceDAE}
\end{equation}

where $E(t)$ is a singular matrix due to rows with only zeros corresponding to the algebraic equations for true sys-LODAEs \footnote{sys-LODEs that contain equations with derivatives in more than one dependent variable may also be written implicitly in this form. In this case, $E(t)$ is not singular.}. This representation is also referred to as semi-state equations, singular systems, or continuous descriptor systems. As is the case with sys-LODEs formulated using state-space, the constitutive laws are not native to this form, and it is not always easy (or direct) to covert them into the state-space formulation. Notice that expressing sys-LODAEs in this particular form leads to dealing with a singular matrix, $E(t)$. This is distinct from, and does not leverage, the fact that the constitutive equations arise from inherently independent physical laws. Consequently, sys-LODAEs formulated using the state-space representation are often referred to as singular systems (especially in the systems and control theory literature) \cite{ZhengBoutat}.

Certain constant coefficient (sys-ccLODAE) are solvable in the time domain or using Laplace transform. The condition for these systems are those for justifying inverting a resultant matrix pencil $sE - A$. For equations of physical systems - which are the systems we are interested in, the matrix pencil is a regular matrix \cite{Milano2017} and hence the problem is solvable \footnote{Generalized matrix inverses (e.g. Drazin inverse) are used for some other resultant matrix pencil.}.

For solution in the time domain, the systems in Equation \ref{eq:stateSpaceDAE} are generally transformed into a canonical form. This reduces the system into two decoupled subsystems in transformed variables.

\begin{equation}
\begin{cases}
\dot{y}(t) = A_1 y(t) + B_1 u(t)\\
N \dot{z}(t) = z(t) + B_2 u(t); \quad N^m = 0, \quad N^{m-1} \neq 0
\end{cases}
\end{equation}

The above can then be solved analytically for sys-ccLODAEs with sufficiently differentiable $u(t)$. This is due to the fact that the first equation is a sys-ccLODE, and the second involves a matrix $N$ that has a special structure and is nilpotent \cite{YipSincovec, ODAEsAreNotODEs}. The nilpotency index is an example of a differential index which is roughly associated with its distance from sys-LODEs and is used to describe the difficulty associated with numerically solving sys-ccLODAEs \cite{NilpotencyIndex}. Specified initial values must be consistent. 

In the $s$-domain, sys-ccLODAE are expressed and solved as generalized eigenvalue problems. Like eigenvalue problems associated with sys-ccLODE, there are existing efforts for developing accurate and efficient algorithms for solving generalized eigenvalue problems which are computationally expensive for large systems \cite{GEPcomputationallyexpensive}. Sys-LODAEs are also studied from a controls perspective for controllability and observability and to design optimal controls  \cite{YipSincovec, Ilchman, ReisVoigt}.

% More generally, sys-ODAEs are typically assigned an index and there are algorithms designed specifically for the IVP numerical simulation of certain lower index. Higher index systems are reduced – using a variety of methods, into lower index systems or sys-ODEs (the conversions are nontrivial and are applicable to specific classes of sys-ODAEs \cite{X}). The resultant sys-LODE IVPs is then solved (often numerically).

Note that existing methods for studying sys-LODEs and sys-LODAEs generally require shifting away from the natural structure of the system. The formulation requires that state variables be defined such that the system can be represented as in Equation \ref{eq:stateSpace1} or \ref{eq:stateSpaceDAE}. %The solutions to state-space problems are generally for initial value problems due to the structure of the state-space formulation and Laplace transform methods. 

%We are not aware of any systematic methods for deriving higher order governing equations from sys-LODAEs.

\clearpageSamiya

\section{Differential Operator Theory for uni-ccLODEs}
\label{s:differentialOperator}

In this section, we review differential operator theory \cite{HungChung} (and introduce modifications) as it is essential in developing our methods for objectives I-III. Classically, differential operator theory methods are extremely powerful and are applied to solve uni-ccLODEs - these are of the form,

\begin{equation}
    P(D) x(t) = f(t)
    \label{eq:appendixUniCCLODE}
\end{equation}

as well as other special uni-ODEs. The complementary solution and particular solution are determined separately. 

Building on its algebraic abstraction and application to solving uni-ccLODEs, we utilize differential operator theory to formulate constitutive equations in Section \ref{s:formulate}, derive governing equations in Section \ref{s:governingMethod} and obtain solutions for a system of sys-ccLODAEs in Section \ref{s:solutionMethod}.

The following operator definitions and rules regarding their manipulation to solve a uni-ccLODE (Equation \ref{eq:appendixUniCCLODE}) are described in the literature \cite{HungChung}. We also make a modification to operator definitions to obtain the full solution more directly, and to more directly handle repeated root cases.

\begin{itemize}
\item \itemEq{Dx(t) \triangleq \frac{d}{dt}x(t)}  
\item \itemEq{Du(t)v(t) = \frac{d}{dt} \big( u(t)v(t) \big) = u(t)Dv(t) + v(t)Du(t) \label{eq:Dto2vars}}
\item \itemEq{P(D)\expn{at} = P(a)\expn{at} \label{eq:toCharacteristicEq}} where $P(.)$ signifies a polynomial in the argument. Classically, this is used to determine the complementary solution of a uni-ccLODE via finding the roots of the characteristic equation $P(a) = 0$.
\item \itemEq{P(D) u(t) v(t)  = u(t) P\bigg(D+\frac{u'(t)}{u(t)}\bigg) v(t)}. Of particular interest is,
\begin{itemize}
    \item \itemEq{D  \expn{at} x(t)  = \expn{at} \big( D + a \big) x(t)}  
\end{itemize}
\item Definition of $\frac{1}{D}$ -  classically and using our modification:
\begin{itemize}
    \item  \itemEq{\frac{1}{D}f(t) \triangleq \defint{t_o}{t}{f(\tau)}{\tau} \label{eq:classicalDefinitionIntegral}}. Is the classical definition in the literature. The lower limit of integration $t_o$ can take on any value and is chosen based on convenience (e.g. $-\infty$). This is because, classically, the constant of integration is neglected. With this, only the particular solution is obtained from Equation \ref{eq:inverseccLODE}.
    \item \itemEq{    \frac{1}{D} f(t) \triangleq \defint{}{}{f(t)}{t} ; \textrm{ with} \quad \frac{1}{D}0=c \label{eq:modifiedDefinitionIntegral}} We redefine Equation \ref{eq:classicalDefinitionIntegral} by mapping $\frac{1}{D}$ onto an indefinite integral and retain constants.  Using this modification, the full solution is obtained from Equation \ref{eq:inverseccLODE} rather than just the particular solution
\footnote{Consider, for instance, the solution to the concocted uni-ccLODE $\big(D^2-3D+2\big) x(t) = \expn{3it}$. Applying Equation \ref{eq:factor} or \ref{eq:partialFractionExpansion} using the classical definition (Equation \ref{eq:classicalDefinitionIntegral}) we obtain only a particular solution, $x(t) = \frac{-7+9i}{130}\expn{3it}$ and must obtain the complementary solution separately. Whereas using the modified definition (Equation \ref{eq:modifiedDefinitionIntegral}), we obtain the full solution $x(t) = c_1 \expn{t} + c_2 \expn{2t} + \frac{-7+9i}{130}\expn{3it}$}. The solution to the homogeneous problem can also be obtained by $x(t) = \frac{1}{P(D)}0$ as an alternative to Equation \ref{eq:toCharacteristicEq}. The modifications also allow for more direct handling of repeated-root cases via Equation \ref{eq:factor} (but not Equation \ref{eq:partialFractionExpansion}) applied to $f(t)=0$ \footnote{For example, solving $(D-a)^2x(t) = f(t)$, with Equation \ref{eq:factor} and the modified definition, very directly yields a complementary solution $(c_1  + c_2 t) \expn{at}$}. 
\end{itemize}

\item \itemEq{\frac{1}{D+a}f(t) = \expn{-at}\frac{1}{D} \expn{at}f(t) \label{eq:l_Da}} which, in turn, is solvable using Definition \ref{eq:classicalDefinitionIntegral} / \ref{eq:modifiedDefinitionIntegral}.
\item \itemEq{x(t) = \frac{1}{P(D)}f(t) \label{eq:inverseccLODE}} is used to solve for the solution of the uni-ccLODE $P(D)x(t) = f(t)$ as follows:
\begin{itemize}
    \item \itemEq{\frac{1}{P(D)}f(t) = \sum\limits_{i=1}^{n} \frac{\gamma_i}{D+\alpha_i} f(t) \label{eq:partialFractionExpansion}} Partial fraction expansion then invoking Equation \ref{eq:l_Da} to obtain the  solution. This expansion is appropriate unless there are repeated roots $\alpha_i$. 
    \item \itemEq{\frac{1}{P(D)}f(t) = \prod\limits_{i=1}^{n} \frac{1}{D+\alpha_i} f(t) \label{eq:factor}} Factorization  and sequentially invoking Equation \ref{eq:l_Da} to obtain the  solution. This is always appropriate. 
\end{itemize}
\item  \itemEq{P(D)\frac{1}{Q(D)} u(t) = \frac{1}{Q(D)}P(D) u(t) ; P(D)Q(D)u(t) = Q(D)P(D)u(t) ;  \frac{1}{P(D)}\frac{1}{Q(D)}u(t) = \frac{1}{Q(D)}\frac{1}{P(D)}u(t) \label{eq:interchangeable}} where $u(t)$ is a smooth function and $P(D), Q(D)$ are polynomials in the differential operator $D$. When applied to smooth functions the order of operations does not matter. This is the case due to the constant nature of polynomial coefficients. 
\end{itemize}

Note that the factorization of $P(D)$ is over complex numbers and hence the irreducible factors are always first order polynomials. To summarize, Equations \ref{eq:partialFractionExpansion} / \ref{eq:factor} are used in conjunction with Equation \ref{eq:l_Da} and Definition \ref{eq:classicalDefinitionIntegral} / \ref{eq:modifiedDefinitionIntegral} to solve uni-ccLODEs $P(D)x(t)=f(t)$. Classically, the complementary and particular solutions are found separately: the complementary solution is obtained from the characteristic equation $P(a)=0$, and the particular solution is obtained $x(t)=\frac{1}{P(D)}f(t)$ by applying \ref{eq:partialFractionExpansion} / \ref{eq:factor} and definition \ref{eq:classicalDefinitionIntegral}. In contrast, using the modified definition (Equation \ref{eq:modifiedDefinitionIntegral}), the \textit{full} solution is obtained from $x(t)=\frac{1}{P(D)}f(t)$ using either Equation \ref{eq:partialFractionExpansion} / \ref{eq:factor}.

% TO DO: XXXXXQ: Has this modification been done before - indefinite rather than definite? Also, handling repeated roots? Any references other than HungChung? Maybe different operator theory exists? Is there such an existing operator theory (abstract or applied) for PDEs? Or is that section new?? See link

%look into TO DO XXXX: https://www.sciencedirect.com/science/article/pii/S074771710290564X

% XXXXX

In what follows, we demonstrate that the modified Definition \ref{eq:modifiedDefinitionIntegral} allows us to arrive at the full solution. We also show that Equation \ref{eq:factor} along with Definition \ref{eq:modifiedDefinitionIntegral} allows us to handle repeated root cases directly. In the main text of this paper, we use Definition \ref{eq:modifiedDefinitionIntegral} (rather than Definition \ref{eq:classicalDefinitionIntegral}), and prefer Equation \ref{eq:factor} to Equation \ref{eq:partialFractionExpansion} for deriving solutions.

\subsection{Full Solution Using Modified Definition}

Using classical methods of differential operator theory (with Definition \ref{eq:classicalDefinitionIntegral}), the complementary and particular solutions are solved for separately and then summed to give the full solution. Using our simple modification (Definition \ref{eq:modifiedDefinitionIntegral}), we may directly obtain the full solution. We do so using Equation \ref{eq:partialFractionExpansion} or \ref{eq:factor}.

We illustrate this by solving the following uni-ccLODE example in detail so that the reader can easily follow and apply our methods for objective III.

\begin{equation}
    \bigg(D^2 + (b-a) D -ab \bigg) x(t) = f(t)
\end{equation}

Solving using Equation \ref{eq:partialFractionExpansion}, we get,

\begin{align}
x(t) &= \frac{1}{(D-a)(D+b)}f(t) \\
 & = \frac{1}{a+b} \bigg( \frac{1}{D-a} - \frac{1}{D+b} \bigg) f(t) \\
 & = \frac{1}{a+b} \bigg( \expn{at} \frac{1}{D} \expn{-at} - \expn{-bt} \frac{1}{D} \expn{bt} \bigg) f(t) \\
 & = \frac{1}{a+b} \bigg( \expn{at} \defint{}{}{\expn{-at} f(t)}{t}  - \expn{-bt} \defint{}{}{\expn{bt}f(t)}{t}  \bigg) 
\end{align}

For instance, for the case of $f(t)=f_o$, this gives:
\begin{align}
   x(t) &= \frac{f_o}{a+b} \bigg( \expn{at} \defint{}{}{\expn{-at} }{t}  - \expn{-bt} \defint{}{}{\expn{bt}}{t}  \bigg) \\
   &= \frac{f_o}{a+b} \Bigg( \expn{at} \bigg(\frac{\expn{-at}}{-a} + c_1\bigg)  - \expn{-bt} \bigg(\frac{\expn{bt}}{b} + c_2 \bigg)  \Bigg) \\
   &= C_1 \expn{at} + C_2 \expn{-bt}
   - \frac{f_o}{ab}
\end{align}

which is the full solution.

\subsection{Handling Repeated Factors}

In some cases, a uni-ccLODE $P(D) x(t) = f(t)$ has a $P(D)$ that, when factorized into first degree polynomial factors, has repeated factors - this corresponds to repeated roots in the characteristic equation of the associated homogeneous equation.

In this situation, the complementary solution (which we can obtain as part of the full solution due to Definition \ref{eq:modifiedDefinitionIntegral}) has a special form. In addition, the uni-ccLODE must be solved using Equation \ref{eq:factor} and cannot be solved using \ref{eq:partialFractionExpansion} (as is obvious by division by zero that is due to the coefficients of partial fraction expansion).

Due to repeated factors, the method of Equation \ref{eq:factor} applies more generally to solving uni-ccLODEs than the method of Equation \ref{eq:partialFractionExpansion} which cannot handle repeated first order factors. Therefore Equation \ref{eq:factor} is our method of choice throughout.

To give show that Definition \ref{eq:modifiedDefinitionIntegral} and Equation \ref{eq:factor} handle repeated factors, consider the following example,
\begin{equation}
    \bigg(D^3 + 2D^2 + D \bigg) x(t) = f(t)
\end{equation}

Again, we solve the uni-ccLODE in detailed steps so that the reader can easily follow and apply our methods for objective III. We get,

\begin{align}
    x(t) &= \frac{1}{D(D+1)(D+1)} f(t) \\
    &= 
    \frac{1}{D}\frac{1}{D+1}\frac{1}{D+1} f(t) \\
    &= 
    \frac{1}{D}\frac{1}{D+1} \expn{-t} \frac{1}{D} \expn{t} f(t)  \\
    &= 
    \frac{1}{D} \expn{-t} \frac{1}{D}  \frac{1}{D} \expn{t} f(t)  \\
    &= \defint{}{}{\expn{-t} \defint{}{}{\defint{}{}{\expn{t} f(t)}{t} }{t} }{t}
\end{align}

For instance, if $f(t)=t$, this gives,
\begin{align}
    x(t) &= \defint{}{}{\expn{-t} \defint{}{}{\defint{}{}{t\expn{t} }{t} }{t} }{t}\\
    &= \defint{}{}{\expn{-t} \defint{}{}{(t \expn{t} - \expn{t} + c_1) }{t} }{t}\\
    &= \defint{}{}{\expn{-t} (c_1 t + t\expn{t} - 2\expn{t} + c_2 )}{t}\\
    &= C_1 \expn{-t} + C_2 t \expn{-t} + C_3 + \frac{t^2}{2} - 2t
\end{align}
where the $t$ in the second term of the complementary solution is due to repeated factors. 

\bibliography{ms}

\end{document}